\documentclass[a4paper,10pt]{article}
\usepackage[utf8]{inputenc}
\usepackage{amsmath}
\usepackage{amssymb}
\usepackage[margin=1.25in]{geometry}
\usepackage{hyperref}\hypersetup{colorlinks,citecolor=blue,linktocpage}

\newcommand{\oh}{\mathcal{O}}

\newcommand{\Z}{\mathbb{Z}}
\newcommand{\R}{\mathbb{R}}
\newcommand{\Q}{\mathbb{Q}}
\newcommand{\G}{\mathbb{G}}

\newcommand{\X}{\mathcal{X}}
\newcommand{\Y}{\mathcal{Y}}
\newcommand{\I}{\mathcal{I}}
\newcommand{\supp}{\operatorname{Supp}}

\newtheorem{theorem}{Theorem}
\newtheorem{definition}{Definition}

\title{On the Singularities of Degenerate Abelian Varieties}
\author{Joseph Tenini~-~University of Georgia~-~jtenini@math.uga.edu}
\date{}
\begin{document}

\maketitle

\begin{abstract}
In 2002, Alexeev provided a modular interpretation for the toroidal compactification
of $A_g$ for the 2nd Voronoi fan. In this paper we show that pairs $(X,\Theta)$ in the boundary of this moduli space are semi-log canonical, the analog of log canonical in the non-normal setting. This extends the result of Koll{\'a}r stating that principally polarized abelian pairs $(X,\Theta)$ are log canonical.
\end{abstract}

\section*{Introduction}

This paper seeks to extend in a natural way the result of Koll{\'a}r \cite[Thm.17.13]{kollar95} which states that principally polarized abelian pairs $(X,\Theta)$ are log canonical. We would like to extend this result to degenerations of such pairs. Though there is more than one compactification of $A_g$, the moduli space of principally polarized abelian varieties, in \cite{alexeev02} Alexeev provides a modular interpretation for the toroidal compactification of $A_g$ for the 2nd Voronoi fan. We will show that the pairs in the boundary of this compactification are in fact semi-log canonical.\\

In section \ref{slc} we give the definition of a semi-log canonical pair. In section \ref{ssp} we review stable semiabelic pairs introduced by Alexeev in \cite{alexeev02}. In section \ref{v} we state a vanishing theorem of Fujino \cite{fujino} which will be used in the proof of the main theorem. In section \ref{r} we state and prove the main result of the paper.\\

\section{Background}

\subsection{Semi-Log Canonical Pairs}\label{slc}

We recall the definitions of log canonical and semi-log canonical from \cite[p.56]{kollarmori98} and \cite[p.193-194]{kollar13}.\\

Let $(X,\Delta)$ be a pair where $X$ is a normal variety, $\Delta = \sum a_iD_i$ is a sum of distinct prime divisors with $a_i \in \Q$, and $K_X + \Delta$ is $\Q$-Cartier. We say that $(X,\Delta)$ is log canonical if the discrepancy of $(X,\Delta)$ is  $\geq -1$. This can be generalized for $X$ non-normal in the following way.\\

\begin{definition}\label{definition}
Let $X$ be a scheme that is $S2$ and whose codimension 1 points are either regular or ordinary nodes (a so called demi-normal scheme). Let $X$ have normalization $\nu:X^{\nu} \to X$ and 
conductors $D \subset X$ and $D^{\nu} \subset X^{\nu}$. Let $\Delta$ be an effective $\Q$-divisor whose support does not contain any irreducible component of $D$ and $\Delta^{\nu}$ the divisorial 
part of $\nu^{-1}(\Delta)$.
The pair $(X,\Delta)$ is called semi-log canonical if $K_X + \Delta$ is $\Q$-Cartier, and 
 $(X^{\nu},D^{\nu}+\Delta^{\nu})$ is log canonical.
\end{definition}

\subsection{Stable Semiabelic Pairs}\label{ssp}

In this paper we consider the principally polarized stable semiabelic pairs introduced by Alexeev in \cite{alexeev02}.
What follows is a summary of the pertinent definitions and results.\\

A semiabelian variety is a group variety $G$ which is an extension
$$ 1 \to T \to G \to A \to 0$$ of an abelian variety $A$ by a torus $T \cong \G_m^r$ (where $\G_m$ is the multiplicative group of the ground field). A semiabelic variety is a variety $P$ 
together with an action of a semiabelian variety $G$ such that:
\begin{enumerate}
 \item $P$ is normal.
 \item The action has only finitely many orbits.
 \item The stabilizer of the generic point is connected, reduced and lies in the toric part $T$ of $G$.
\end{enumerate}

In the special case that $G$ is abelian, we say that $P$ is an abelic variety.\\

A stable semiabelic variety is a variety $P$ along with an action of 
a semiabelian variety $G$ of the same dimension satisfying the same conditions except condition (1) is weakened to require only seminormality. We say that a reduced scheme $P$ is seminormal if every proper bijective morphism $f:P' \to P$ with reduced $P'$ inducing isomorphisms on the residue fields $\kappa(p') \supset \kappa(p)$ for 
each $p \in P$ is an isomorphism.\\

A stable semiabelic variety is polarized if it is projective and comes with an ample invertible sheaf $L$. The degree of the polarization is $h^0(L)$.
A stable semiabelic pair $(P,\Theta)$ is a stable projective semiabelic variety $P$ and an effective ample Cartier divisor $\Theta$ which does not contain any $G$-orbits. $P$ is polarized by 
$L = \oh_P(\Theta)$. We define an abelic pair analogously as a projective abelic variety (also called an abelian torsor) with an effective ample Cartier divisor.\\

In \cite{alexeev02} we see that such pairs correspond to certain complexes of lattice polytopes referenced by $\overline{\Lambda}/Y$ where $\Lambda \cong \Z^g$ is a lattice, $Y \subset \Lambda$ is 
a subgroup and $\overline{\Lambda}$ is a trivial $\Lambda$-torsor. Analysis of such pairs is categorized by type which is determined by $\overline{\Lambda}$ and $Y$. For this paper we will be 
concerned with the principally polarized case, which corresponds to the case when $|\Lambda/Y| = 1$.\\

Alexeev establishes the following facts concerning the moduli of abelic and semiabelic pairs:

\begin{theorem}\label{alexeevthm1}
\begin{enumerate}
\item The moduli stack $\mathcal{AP}_{g,d}$ of abelic pairs $(P,\Theta)$ of degree $d$ is a separated Artin stack with finite stabilizers and comes with a natural map of relative dimension $d-1$ 
to the stack $\mathcal{A}_{g,d}$ of polarized abelian varieties \cite[Thm.1.2.2]{alexeev02}.
\item $\mathcal{AP}_{g,d}$ has a coarse moduli space $AP_{g,d}$ which is a separated scheme and comes with a  natural projective map of relative dimension $d-1$ to $A_{g,d}$ \cite[Thm.1.2.2]{alexeev02}.
\item In the principally polarized case (that is, when $d = 1$), we have $\mathcal{AP}_g = \mathcal{A}_g$ and $AP_g = A_g$ \cite[Thm.1.2.3]{alexeev02}.
\item The component $\overline{\mathcal{AP}}_{g,d}$ of the moduli stack of semiabelic pairs containing $\mathcal{AP}_{g,d}$ and pairs of the same numerical type is a proper Artin stack with finite stabilizers \cite[Thm.1.2.16]{alexeev02}.
\item It has a coarse moduli space $\overline{AP}_{g,d}$ as a proper algebraic space \cite[Thm.1.2.16]{alexeev02}.
\item In the principally polarized case (when $Y=\Lambda$) the toroidal compactification of $A_g$ for the second Voronoi decomposition is isomorphic to the normalization of the main irreducible 
component of $\overline{AP}_{g,1}$, the one containing $A_g = AP_g$ \cite[Thm.1.2.17]{alexeev02}.
\end{enumerate}
\end{theorem}

In the above paper there are also many facts about principally polarized stable semiabelic pairs, some of which are implicit in the arguments. Below we summarize and make explicit the ones that we will use in our analysis.

\begin{theorem}\label{alexeevthm2}
Let $(\X_0,\Theta_0)$ be a stable semi-abelic pair that appears as the special fiber of a flat family $\pi:(\X,\Theta) \to S$ where $S$ is a smooth curve and a general fiber of $\pi$ is a principally polarized abelian variety.
\begin{enumerate}
\item $K_{X_0} \sim 0$, where $\sim$ denotes linear equivalence.
\item There is an $\epsilon > 0$ such that $(X_0,\epsilon\Theta_0)$ is semi-log canonical.
\item $h^0(X_0,\oh_{X_0}(\Theta_0)) = |\Lambda/Y| = 1$.
\item $H^i(X_0,\oh_{X_0}(\Theta_0))=0$ for all $i > 0$ \cite[Thm.5.4.1]{alexeev02}.
\end{enumerate}
\end{theorem}

These statements require some remarks.\\

(1) By the work of \cite{mumford72}, there is a toric model of $(X,X_0)$. That is, we may write $(X,X_0)=(Y,Y_0)/\Z^r$ where $Y$ is toric and $Y_0$ is the toric boundary. Thus, by adjunction $K_{Y_0} \sim 0$ and so, since the action is $\Z^r$-equivariant, we also have $K_{X_0} \sim 0$.\\

(2) By \cite[Thm.1.2.14]{alexeev02}, $X_0$ is Cohen-Macaulay and hence S2 and so the existence of an $\epsilon > 0$ such that $(X_0,\epsilon\Theta_0)$ is semi-log canonical amounts to showing that (in the notation of definition \ref{definition}) $K_{X_0}+\Theta_0$ is $\Q$-Cartier, $(X_0^{\nu},D^{\nu})$ is log canonical and $\Theta_0^{\nu}$ does not contain any log canonical centers of $(X_0^{\nu},D^{\nu})$. $K_{X_0}+\Theta_0$ is $\Q$-Cartier since $\Theta_0$ is Cartier and $K_{X_0} \sim 0$. $(X_0^{\nu},D^{\nu})$ is log canonical because $K_{X_0^{\nu}}+D^{\nu} \sim 0$ and $(X_0^{\nu},D^{\nu})$ is a toric pair.
Finally, by definition, $\Theta_0$ does not contain any $G$-orbits, and the log canonical centers of $(X_0^{\nu},D^{\nu})$ are precisely the closures of the codimension $\geq 1$ orbits. This implies in particular that $X_0$ is a demi-normal scheme such that $\Theta_0$ does not contain any irreducible component of the conductor of normalization. By continuity of discrepancies, the pair $(X_0,\epsilon \Theta_0)$ is semi-log canonical for $0 < \epsilon \ll 1$.\\

(3) Since we work with principally polarized pairs, the condition $h^0(X,\oh_X(\Theta)) = |\Lambda/Y| = 1$ is satisfied by definition.\\

(4) This vanishing result is the most important property that we will need in our proof of the main result. A complete proof of this fact is given in \cite[Thm.5.4.1]{alexeev02}.\\

\subsection{Vanishing}\label{v}

In Koll$\operatorname{\acute{a}}$r's proof for principally polarized abelian varieties \cite{kollar95}, he uses the vanishing theorem of Kawamata and Viehweg. We will proceed analogously but using a vanishing theorem of Fujino \cite[Thm.6.3]{fujino}, which we use only with $\Q$-divisors.\\

\begin{theorem}\label{vanish}
Let $Y$ be a smooth variety and let $B$ be a boundary $\R$-divisor such that $\supp B$ is simple normal crossing. Let $f: Y \to X$ be a projective morphism and let $L$ be a Cartier divisor on $Y$ such that $L - (K_Y + B)$ is $f-semi-ample$. Let $\pi: X \to S$ be a projective morphism. Assume that $L - (K_X+B) \sim_{\R} f^* H$ for some $\pi$-ample $\R$-Cartier $\R$-divisor $H$ on $X$. Then $R^p\pi_*R^qf_*\oh_Y(L)=0$ for every $p > 0$ and $q \geq 0$.
\end{theorem}

We are now ready to state and prove the main result of the paper.

\section{Statement and Proof of Main Result}\label{r}

We work now over an algebraically closed field of characteristic zero.

\begin{theorem}
Let $(\X_0,\Theta_0)$ be a stable semi-abelic pair that is a one-parameter degeneration. That is, we have a flat family $\pi:(\X,\Theta) \to S$ where $S$ is a smooth curve, a general fiber of $\pi$ is a principally polarized abelian variety, and the special fiber of $\pi$ is $(\X_0,\Theta_0)$. Then $(\X_0,\Theta_0)$ is semi-log canonical.
\end{theorem}

\textbf{Proof:} As we have seen in Theorem \ref{alexeevthm2} we already have that $\X_0$ is demi-normal, $K_{\X_0}+\Theta_0$ is Cartier, and $\Theta_0$ does not contain any irreducible component of the conductor of normalization. Thus, by adjunction $(\X_0,\Theta_0)$ is semi-log canonical if $(\X,\Theta + \X_0)$ is log canonical in a neighborhood of $\X_0$. So, it would suffice to 
show that $(\X, \Theta + (1-\epsilon)\X_0)$ is log canonical for every $0 < \epsilon \ll 1$.\\

Let $f: \Y \to \X$ be a resolution of singularities so that $\Y$ is smooth and $\widehat{\X_0} \cup \widehat{\Theta} \cup E$ is a simple normal crossing divisor where 
$\widehat{\X_0} = f^{-1}_*(\X_0)$, $\widehat{\Theta} = f^{-1}_*(\Theta)$, and $E$ is the exceptional locus of $f$. Such a resolution is guaranteed by Hironaka \cite{hironaka}.\\

Fix $0 < \epsilon \ll 1$ and put $d = 1 - \epsilon$. Let $c$ be the log canonical threshold so that $(\X, c(\Theta + d\X_0))$ is log canonical but $(\X, c'(\Theta + d\X_0))$ is not log canonical for any $c' > c$. 

Seeking a contradiction, suppose that $c < 1$ and put $cd = d'$. This means that we may write the following equation:

$$ K_{\Y} + c\widehat{\Theta} + d'\widehat{\X_0} = f^*(K_{\X} + c\Theta + d'\X_0) + D $$

\noindent where $D$ is exceptional. Write $D$ as $D = B - A - \Delta$ where $A,B,\Delta$ are effective, $A,B$ are integral and nonzero, $\lfloor \Delta \rfloor = 0$, and $A,B$ have no irreducible components in common.

Now the equation becomes:

$$ K_{\Y} + c\widehat{\Theta} + d'\widehat{\X_0} = f^*(K_{\X} + c\Theta + d'\X_0) + B - A - \Delta $$

\noindent Which we can rewrite as:

$$ -f^*(K_{\X} + c\Theta + d'\X_0) = B - A - (K_{\Y} + c\widehat{\Theta} + d'\widehat{\X_0} + \Delta) $$

\noindent Adding $f^*(K_{\X} + \Theta + K_{\X_0})$ to each side:

$$ f^*((1-c)\Theta + (1-d')\X_0) = (f^*(K_{\X} + \Theta + \X_0) + B - A) - (K_{\Y} + c\widehat{\Theta} + d'\widehat{\X_0} + \Delta) $$

\noindent Now we would like to apply Theorem \ref{vanish} to the above situation. The conclusion is: 
\begin{equation}\label{eq1}\tag{$\star$}
 R^p\pi_* R^q f_* \oh_{\Y}(f^*(K_{\X} + \Theta + \X_0) + B - A) = 0 
\end{equation}
\noindent for $ q \geq 0$ and $p > 0$.\\

Since $B$ is exceptional and $A$ is non-zero, $f_*\oh_{\Y}(B - A) \subset \oh_{\X}$ is a non-trivial ideal sheaf. Call this ideal $\I$ and let $Z \subset \X$ be the subscheme determined by $\I$.\\

Consider the short exact sequence:

$$ 0 \to \I \to \oh_{\X} \to \oh_Z \to 0 $$

\noindent Tensoring by $\oh_{\X}(\Theta)$ we get:

$$ 0 \to \oh_{\X}(\Theta) \otimes \I \to \oh_{\X}(\Theta) \to \oh_{Z}(\Theta) \to 0 $$

\noindent Applying $\pi_*$, we arrive at the long exact sequence:

$$ 0 \to \pi_*(\oh_{\X}(\Theta) \otimes \I) \to \pi_*(\oh_{\X}(\Theta)) \to \pi_*(\oh_{Z}(\Theta)) \xrightarrow{\rho} R^1\pi_*(\oh_{\X}(\Theta) \otimes \I) \to \dots$$

\noindent By \eqref{eq1}, we have $R^1\pi_*(\oh_{\X}(K_{\X} + \Theta + \X_0) \otimes \I) = 0$ and so $R^1\pi_*(\oh_{\X}(\Theta) \otimes \I)=0$. Thus, we have that $\rho: \pi_*(\oh_{\X}(\Theta)) \to \pi_*(\oh_{Z}(\Theta))$ is surjective.\\

Let us now consider the sheaf $\pi_*\oh_{\X}(\Theta)$. From Theorem \ref{alexeevthm2} Part (4), we have that $H^i(\X_0,\oh_{\X_0}(\Theta))=0$ for all $i > 0$. Thus, by the cohomology and base change formula (see for example \cite[III.12.11]{hartshorne}) $\pi_*\oh_{\X}(\Theta)$ is locally free of rank 1. Moreover, all sections of $\pi_*(\oh_{\X}(\Theta))$ vanish on $Z$ (since $Z \subset \Theta$). Thus, $\rho$ is in fact the zero map since it is merely the restriction to $Z$. We will reach a contradiction if we can show that $\pi_*(\oh_{Z}(\Theta))$ is not the zero sheaf.\\

Indeed, if $\tau$ is a general translation of $\X$, we will have $\tau(\Theta) \not\subset Z$. Specifically, $\tau$  is a section of a semiabelian scheme $G$ over $S$. That is, $G$ is a smooth separated group scheme over $S$ such that every geometric fiber $G_{p}$ is a semiabelian variety over $\kappa(p)$. There is a question as to the existence of such a section, and indeed, such a section may not exist. However, a section does exist after an extension $G \times_S S' \to S'$. Taking $\tau$ to be such a section we have $\tau(\Theta) \not\subset Z$ as desired.
Hence, $H^0(\X,\pi_*(\oh_{Z}(\tau\Theta)))\not= 0$ and so by semi-continuity $H^0(\X,\pi_*(\oh_{Z}(\Theta)))\not= 0$ as well. We have thus arrived at a contradiction and the proof is complete. \begin{flushright}$\Box$ \end{flushright}

\section*{Acknowledgments} I would like to thank Valery Alexeev for raising the question that is the subject of this paper as well as providing many insightful comments and suggestions.

\bibliography{TeniniDegenerations}
\bibliographystyle{plain}

\end{document}